\documentclass[11pt]{article}
\usepackage{latexsym}
\usepackage{amsfonts}
\usepackage{amsmath}
\makeatletter

\def\@footnotetext#1{\insert\footins{%

\footnotesize 

 \interlinepenalty\interfootnotelinepenalty

 \splittopskip\footnotesep

 \splitmaxdepth \dp\strutbox \floatingpenalty \@MM

 \hsize\columnwidth \@parboxrestore

 \edef\@currentlabel{\csname 
p@footnote\endcsname\@thefnmark}\@makefntext
 {\rule{\z@}{\footnotesep}\ignorespaces
#1\strut}}}

\def\abstract{\small\quotation{\hskip-\parindent\sc Abstract.}}
\def\classification{\@ifnextchar [{\@xfootnotenext}%
 {\begingroup\let\protect\noexpand
 \xdef\@thefnmark{}\endgroup
 \@footnotetext}}

\title {}

\begin{document}

\begin{center}

\classification {{\it 2000 Mathematics Subject Classification:} 
Primary 14R10, 14R15; Secondary 14A10, 13B25.\\
$\ast$) Partially supported by RGC Grant Project HKU 7134/00P.}

{\bf \Large  Test polynomials, retracts, and  the Jacobian conjecture}

\bigskip

{\bf Vladimir Shpilrain }

\medskip

 and

\medskip

 {\bf Jie-Tai Yu}$^{\ast}$

\end{center}

\medskip

\begin{abstract}
\noindent  Let $K[x,y]$ be the algebra of 
two-variable polynomials over a field $K$. A polynomial $p=p(x, y)$ is 
called a {\it test polynomial} 
for automorphisms if, whenever $\varphi(p)=p$ for a mapping 
 $\varphi$ of $K[x,y]$, this  $\varphi$ must be an automorphism. 
Here we show that $p \in {\mathbb C}[x,y]$ is 
a  test polynomial if and only if $p$  does not belong to any 
proper retract of ${\mathbb C}[x,y]$. 
This  has the following  corollary that may have 
application to the Jacobian conjecture: if a  mapping  $\varphi$ of
 ${\mathbb C}[x,y]$ with invertible Jacobian matrix is ``invertible on 
one particular polynomial", then 
it is an automorphism. More formally: if there is a non-constant polynomial $p$ and 
 an injective mapping  $\psi$ of ${\mathbb C}[x,y]$  such that 
$\psi(\varphi(p))=p$, then $\varphi$ is an automorphism. 

\end{abstract} 

\bigskip

\section{Introduction}

Let $K[x,y]$ be the algebra of 
two-variable polynomials over a field $K$ of characteristic 0. 
 A subalgebra $R$ of $K[x,y]$ is called a {\it retract} if there 
is an
idempotent homomorphism $\pi$ of $K[x,y]$
(called a {\it retraction} or a {\it projection})
such that $\pi(K[x,y])=R$. 

 There are several equivalent descriptions of retracts of $K[x,y]$ 
known by now:
\smallskip

{\bf (i)} $K[x,y]=R\oplus I$ for some ideal $I$ of $K[x,y]$;
\smallskip

{\bf (ii)} $K[x,y]$ is a projective extension of $R$ in the category 
of $K$-algebras;
\smallskip

{\bf (iii)} By a theorem of Costa \cite{Costa}, every proper retract 
of $K[x,y]$ 
(i.e., one different from $K[x,y]$ and $K$) is of the form 
$K[p]$ for some $p=p(x,y)\in K[x,y]$.  The authors earlier 
proved \cite{ShYu} that there exists an automorphism of $K[x,y]$ which 
takes $p(x,y)$ to   $x+y\cdot q(x,y)$ for some $q(x,y)\in K[x,y]$, and
every polynomial of the form $x+y\cdot q(x,y)$
generates a proper retract of $K[x,y]$.  
\smallskip

{\bf (iv)} (see \cite{ShYu})   
$p(x,y)$ generates a retract of $K[x,y]$
if and only if there is an endomorphism of $K[x,y]$ which takes
$p(x,y)$ to $x$.
\smallskip

{\bf (v)} (see \cite{DrenskyYu}) $p(x,y)$ belongs to a proper 
retract of ${\mathbb C}[x,y]$ if and only if $p(x,y)$ is fixed by some 
endomorphism of ${\mathbb C}[x,y]$ with nontrivial kernel.

\smallskip

Recently retracts have found another application in a general setup 
of arbitrary free algebras and groups in relation with test elements, 
introduced in \cite{Sh}. In general, an element $g$ of a group or
an algebra $F$ is a  test element if any endomorphism of $F$ fixing 
$g$ is actually an automorphism. It is easy to see that a test element  
does not belong to any proper retract of $F$; a remarkable result of 
Turner  \cite{Turner} says that, if $F$ is a free group, then the 
converse is also true. Thus, an element of a free group $F$ 
is a  test element if and only it  does not belong to any 
proper retract of $F$. 

 Here we establish a similar characterization of test polynomials 
in ${\mathbb C}[x,y]$:

\medskip 

\noindent {\bf Theorem  1.} A polynomial $p \in {\mathbb C}[x,y]$ is 
a  test polynomial if and only if $p$  does not belong to any 
proper retract of ${\mathbb C}[x,y]$.

\medskip 

 Our proof uses several recent results, in particular, a result of 
Drensky and  Yu \cite{DrenskyYu} mentioned in the item (v) above. 
Crucial for 
our proof is the following result of independent interest. 

\medskip 

\noindent {\bf Theorem  2.}  Let $\varphi$ be an injective 
endomorphism of ${\mathbb C}[x,y]$ which is not an automorphism. 
Suppose that $\varphi(p)=p$ for some non-constant polynomial 
$p \in {\mathbb C}[x,y]$. 
Then $p \in {\mathbb C}[q]$, where $q$ is a coordinate polynomial 
of ${\mathbb C}[x,y]$. In particular, $p$ belongs to a proper retract of ${\mathbb C}[x,y]$.
\medskip 

 Recall that $q=q(x, y)$ is a {\it coordinate polynomial} of ${\mathbb C}[x,y]$ 
if it can be taken to $x$ by an automorphism of ${\mathbb C}[x,y]$. 

 We also use results of Shestakov and 
Umirbaev \cite{Umirbaev} on estimating degrees of polynomials in 
two-generated subalgebras of $K[x,y]$.   Another   ingredient is 
a result of Kraft \cite{Kraft} concerning the subalgebra
 $\varphi^\infty({\mathbb C}[x,y]) = 
\cap_{k=1}^\infty \varphi^\infty({\mathbb C}[x,y])$.

Theorems  1, 2 have the following  corollary:
\medskip 

\noindent {\bf Corollary.} Let $\varphi$ be an  
endomorphism of ${\mathbb C}[x,y]$ with invertible Jacobian matrix. 
 If there is a non-constant polynomial $p \in {\mathbb C}[x,y]$ and 
 an injective mapping  $\psi$ such that $\psi(\varphi(p))=p$, then 
$\varphi$ is an automorphism of  ${\mathbb C}[x,y]$.
\medskip 

 This   strengthens our earlier result \cite[Corollary 1.7]{ShYu}, where 
we showed that, if $\varphi$ has invertible Jacobian matrix, then 
$\varphi(p)=p$ implies that $\varphi$ is an automorphism of  ${\mathbb C}[x,y]$.
\smallskip

 To conclude the Introduction, we raise a problem motivated by 
results of this paper:

\medskip 

\noindent {\bf Problem.} Suppose 
$p \in {\mathbb C}[x,y]$ is a test polynomial and 
$\varphi$ is an injective mapping of ${\mathbb C}[x,y]$. Is 
 $\varphi(p)$ necessarily a test polynomial?  
\medskip 

 It is interesting to note that, by a result of Jelonek \cite{Jelonek}, 
a ``generic"  
polynomial of degree $\geq 4$ is a test polynomial.

\section{Proof of Theorem  2}

 We consider the following two principal cases.

\smallskip

\noindent {\bf  Case I.} There is a coordinate polynomial in 
$\varphi({\mathbb C}[x,y])$. 
\smallskip

\noindent {\bf  Case II.} There are no coordinate polynomials in 
$\varphi({\mathbb C}[x,y])$. 
\smallskip

 In Case I, consider two subcases: 
\smallskip

\noindent {\bf (1)} $\varphi$ is {\it not birational}, i.e., 
does not induce an automorphism of the field of fractions. 

Then, by a result of Kraft \cite[Lemma 1.3]{Kraft}, 
$\varphi^\infty({\mathbb C}[x,y])$ is 
either ${\mathbb C}$ or ${\mathbb C}[f]$,  where $f=f(x, y)$ is
some polynomial. Obviously, if $\varphi(p)=p$, then 
$p \in \varphi^\infty({\mathbb C}[x,y])$. 
We are therefore going to focus on the case 
$\varphi^\infty({\mathbb C}[x,y]) = {\mathbb C}[f]$ and show that, if 
$\varphi$ is injective, then $p \in {\mathbb C}[q]$, where $q$ is 
a coordinate polynomial. 

 Now suppose $r=r(x,y)$ is a coordinate polynomial in 
$\varphi({\mathbb C}[x,y])$, and let $r=\varphi(s(x,y))$. Then, 
since $\varphi$ is injective, the polynomial $s=s(x,y)$ must be 
coordinate, too, by the result of  \cite{Campbell}. 
Therefore, upon changing generating set of ${\mathbb C}[x,y]$ if necessary, 
we may assume that $r=\varphi(x)$. Furthermore, we can  
replace $\varphi$ with its conjugate by 
an arbitrary automorphism, say $\alpha$, i.e., with 
$\psi = \alpha \varphi \alpha^{-1}$, and at the same time 
replace $p$ with $p_1=\alpha(p)$. Then we have:

 $$\psi(p_1) = \alpha \varphi \alpha^{-1}(\alpha(p)) = \alpha(p) =  p_1.$$

 Therefore, the pair $(\psi, p_1)$ has the same properties that 
 the pair $(\varphi, p)$ does, namely, $\psi$ is injective but 
not birational,  and $\psi(p_1) = p_1$; in particular,  
$p_1 \in \psi^\infty({\mathbb C}[x,y])$. 
By  choosing  $\alpha$ appropriately, we can also have  $\psi(x)=x$, 
thus getting $x \in \psi^\infty({\mathbb C}[x,y])$. Then,  if 
$p$ (and therefore $p_1$) does not belong to ${\mathbb C}[q]$ for any
coordinate polynomial $q$, $\psi^\infty({\mathbb C}[x,y])$ cannot be 
of the form ${\mathbb C}[f]$, which is in  contradiction with the 
result of Kraft mentioned above. 
This completes case (1).

\smallskip

\noindent {\bf (2)} $\varphi$ is  birational, i.e., 
  induces an automorphism of the field of fractions. Again, as in the case (1) 
above, we deduce from \cite{Campbell} that $\varphi$ must take some 
coordinate polynomial to coordinate. Thus, 
 upon changing generating set of ${\mathbb C}[x,y]$ if necessary, 
we may assume that $\varphi$ takes $x$ to $u$, and $y$ to $v \cdot f(u)$, 
where ${\mathbb C}[u,v] = {\mathbb C}[x,y]$, and $f(u)$ is a non-constant polynomial 
(otherwise, $\varphi$ would be an automorphism). 

 Now let 
$$p = p(x,y) = \sum_{i, j} c_{ij} x^i y^j. $$

 Then 
$$\varphi(p) = \sum_{i, j} c_{ij} u^i v^j (f(u))^j. $$

 Let $x^r y^s$ be the highest term of $p(x,y)$ in the ``pure lex" order 
with $y > x$. Then in $\varphi(p)$, the highest term is that of $\varphi(x^r y^s)$
because the $y$-degree of $\varphi(y)$ is not lower than that of $\varphi(x)$.
Furthermore, the highest term of $\varphi(x^r y^s)$ must have the $y$-degree 
at least $s$ since otherwise, one would have both $u$ and $v$ of 
$y$-degree equal to 0, which is impossible. 

 If the $y$-degree of $\varphi(x^r y^s)$ is $>s$, this gives a contradiction 
with $\varphi(p) = p$. 
 Now suppose the $y$-degree of $\varphi(x^r y^s)$ is {\it exactly} $s$. 
This is only possible if the $y$-degree of $v$ is 1 and the $y$-degree of $u$ is 0.
Then, arguing as in the case (1) above, we may assume that $\varphi(x) = x$. 
Therefore, $\varphi(y) = (y+g(x)) \cdot f(x)$. Then from $\varphi(p) = p$ we get:

 $$\sum_{i, j} c_{ij} x^i y^j =  \sum_{i, j} c_{ij} x^i (y+g(x))^j (f(x))^j.$$

 Again we use the ``pure lex" order with $y > x$ to focus on the monomial of 
highest degree on either side, but this time we compare the $x$-degrees of these 
highest-degree monomials. We see that these $x$-degrees cannot be equal unless 
$f(x)$ is a constant, contradicting the assumption. This completes the proof 
in Case I. 

\smallskip

 In Case II, we are going to prove the following somewhat stronger 
statement:

\medskip 

\noindent {\bf Proposition.}  Let $\varphi$ be an injective 
endomorphism of ${\mathbb C}[x,y]$, and 
suppose that there are no coordinate polynomials in 
$\varphi({\mathbb C}[x,y])$. Then $\varphi^\infty({\mathbb C}[x,y]) = {\mathbb C}$. 

\medskip 

\noindent {\bf Proof.}  Let  $\varphi(x) = u = u(x,y), 
~\varphi(y) = v = v(x,y)$, and let $D(u, v)$ denote the determinant 
of the Jacobian matrix of $\varphi$. Since $\varphi$ is injective, 
$D(u, v) \ne 0$. Now there are two cases:
\smallskip

\noindent {\bf  (1)} deg($D(u, v)$)$=0$, i.e., $D(u, v)$ is a 
 non-zero constant. Then, by a result of Kraft \cite{Kraft}, we have 
$\varphi^\infty({\mathbb C}[x,y]) = {\mathbb C}$. 
\smallskip

\noindent {\bf  (2)} deg($D(u, v)$)$>0$.  
Note that for any $k \ge 1$, there are no coordinate polynomials in 
$\varphi^k({\mathbb C}[x,y])$. 
Indeed, if there were a coordinate polynomial in $\varphi^k({\mathbb C}[x,y])$, 
then, by the result of \cite{Campbell}, there would have to be a 
coordinate polynomial in $\varphi^{k-1}({\mathbb C}[x,y])$. 
This would lead to a contradiction with the assumption that there 
are no coordinate polynomials in $\varphi({\mathbb C}[x,y])$.
 
Let  $\varphi^k(x)=u^{(k)}, ~\varphi^k(y) = v^{(k)}$. Then from deg($D(u, v)$)$>0$ 
and from the ``chain rule" we get deg($D(u^{(k)}, v^{(k)})$) $\ge k$. 
Now the Proposition will follow from the lemma below. Before we get 
to it, we need one more definition. 

We call a pair $(p, q)$ of polynomials from $K[x,y]$ {\it elementary reduced} if 
the sum of their  degrees cannot be reduced by a  (non-degenerate) linear transformation
 or a transformation of one of  the following two types: 
\smallskip

\noindent {\bf (i)} $(p, q) \longrightarrow (p+ \mu \cdot 
q^k, q)$ for some $\mu \in K^\ast; ~k \ge 2$;
\smallskip

\noindent {\bf (ii)}  $(p, q) \longrightarrow 
(p, q+ \mu \cdot p^k)$. 
\smallskip

 Now we are ready for our 

\medskip 

\noindent {\bf Lemma.} Let  $p=p(x,y)$ and $q=q(x,y)$ be two 
algebraically independent polynomials such that the pair 
$(p, q)$ is elementary reduced. Let $n=\deg(p) <  m=\deg(q)$; 
$m, n \ge 2$, deg($D(p, q)$) $\ge k$. 
Let  $w=w(x,y) \in {\mathbb C}[p,q])$. Then, unless $w$ is a linear combination 
of $p$ and $q$, one has $\deg(w) > min(n, k).$
\medskip 

\noindent {\bf Proof.} The proof here is based on a result of 
Shestakov and  Umirbaev \cite[Theorem 3]{Umirbaev}.  
Let $N=N(p,q)=\frac{mn}{g.c.d.(n,m)}-m-n+\deg(D(p, q))+2$. 
Following \cite{Umirbaev}, we may assume that the highest homogeneous parts 
of $p$ and $q$ are algebraically dependent; otherwise,  $\deg(w) > n$ 
is immediate (unless $w$ is a linear combination of $p$ and $q$). 
Then $\frac{mn}{g.c.d.(n,m)}-m-n \ge 0.$ Indeed, if $g.c.d.(n,m)=n$, then 
the pair $(p, q)$ would not be elementary reduced, contradicting the assumption.
If $g.c.d.(n,m)<n$, then $\frac{n}{g.c.d.(n,m)} \ge 2$, therefore 
$\frac{mn}{g.c.d.(n,m)} \ge 2m$, hence $\frac{mn}{g.c.d.(n,m)}-m-n \ge 0.$

 Thus, from now on we assume $N=N(p,q) \ge \deg(D(p, q))+2$. 

Suppose now that the $y$-degree of $w=w(x,y)$ is of the form 
$\frac{n}{g.c.d.(n,m)} \cdot b + r \ne 0$, where $0 \le r < \frac{n}{g.c.d.(n,m)}$. 
 Then, by \cite[Theorem 3]{Umirbaev}, we have 

$$\deg(w(p, q)) \ge b \cdot N + mr.$$ 

 If $b \ne 0$, this implies $\deg(w(p, q)) \ge N \ge k+2 >k.$ 
If $b = 0$, then $r \ne 0$, implying $\deg(w(p, q)) \ge m > n.$

 It remains to consider the case where the $y$-degree of $w=w(x,y)$ is 0. 
Then the $x$-degree of $w$ must be nonzero; suppose it is of the form 
$\frac{m}{g.c.d.(n,m)} \cdot b_1 + r_1 \ne 0$, where $0 \le r_1 < \frac{m}{g.c.d.(n,m)}$. 
 Then, again by \cite[Theorem 3]{Umirbaev}, we have 

$$\deg(w(p, q)) \ge b_1 \cdot N + nr_1.$$ 

 As before, $b_1 \ne 0$ implies $\deg(w(p, q))  >k.$ If $b_1 = 0$, then $r \ge 2$ 
because we assume that $w(x,y)$ is not linear. Then we have $\deg(w(p, q)) \ge 2n > n$, 
which completes the proof of the lemma. $\Box$ 

\smallskip

 Continuing with the proof of the Proposition, we aim at showing that  
for any integer $M$, there is an integer $k$ such that the degree of 
any polynomial in $\varphi^k({\mathbb C}[x,y])$ is $>M$. 
The above lemma ``almost" does it if we use it with $p=\varphi^k(x)=u^{(k)}, 
~q=\varphi^k(y) = v^{(k)}$, but it has one extra condition 
on the pair $(p, q)$ to be elementary reduced, whereas a pair 
$(u^{(k)}, v^{(k)})$ may not be elementary reduced. However, if we 
denote by $({\overline u^{(k)}}, {\overline v^{(k)}})$ an elementary reduced 
pair obtained from $(u^{(k)}, v^{(k)})$ by elementary  transformations, 
we shall have all conditions of the lemma satisfied for this pair while 
obviously ${\mathbb C}[{\overline u^{(k)}}, {\overline v^{(k)}}] = {\mathbb C}[u^{(k)}, v^{(k)}].$ 
 In particular, the inequality deg($D({\overline u^{(k)}}, {\overline v^{(k)}})$) $\ge k$ 
 follows from the fact that the mapping $x \to {\overline u^{(k)}}, y \to {\overline v^{(k)}}$ 
is a composition of $\varphi^k$ with an automorphism $\alpha$ of ${\mathbb C}[x,y]$ 
in such a way that $\alpha$ is applied first. Therefore, the ``chain rule" 
applied to this  composition yields deg($D({\overline u^{(k)}}, {\overline v^{(k)}})$)=
deg($D(u^{(k)}, v^{(k)})$). Thus, our lemma is applicable to the pair 
$({\overline u^{(k)}}, {\overline v^{(k)}})$, which completes the proof of the 
Proposition and therefore of Theorem  2. $\Box$

\section{Proof of Theorem 1 and  Corollary}

 The ``only if" part of Theorem 1 follows from a result of \cite{ShYu} 
rather easily. If $p=p(x,y)$ belongs to a proper retract ${\mathbb C}[q]$ of 
${\mathbb C}[x,y]$, then, by \cite{ShYu}, for some automorphism $\alpha$, 
$\alpha(p)$ belongs to ${\mathbb C}[x+y\cdot u]$ for some polynomial $u=u(x,y)$. 
Then the mapping $x \to x+y\cdot u, ~y \to 0$ fixes the polynomial $x+y\cdot u$, 
and therefore also fixes $\alpha(p)$. Thus, $\alpha(p)$ is not a test polynomial, 
and  neither is $p$.

 For the ``if" part of Theorem 1, suppose that $p$ does not belong to 
any proper retract of ${\mathbb C}[x,y]$, and let $\varphi(p)=p$ for some 
mapping $\varphi$ of ${\mathbb C}[x,y]$. Then, by the result of \cite{DrenskyYu}, 
$\varphi$ must be injective. Then, by our Theorem 2, $\varphi$ must be 
an automorphism, hence $p$ is a test polynomial. $\Box$
\medskip 

\noindent {\bf Proof of Corollary 1.} By way of contradiction, assume that 
$\varphi$ is not an automorphism. Then, by our Theorem 2, $p \in {\mathbb C}[q]$, 
where $q$ is a coordinate polynomial of ${\mathbb C}[x,y]$. Therefore, 
the composite mapping $\psi \varphi$ fixes a polynomial $f(q)$ in $q$. 
Then it is easy to see (by looking at the highest degree monomial in $f(q)$) 
that $\psi(\varphi(q))= c \cdot q$ for some $c \in {\mathbb C}^\ast$, 
which implies, by the result of \cite{Campbell}, that $\varphi(q)$ is 
coordinate. A mapping of ${\mathbb C}[x,y]$ with invertible Jacobian matrix 
that takes a coordinate polynomial to a coordinate polynomial is obviously 
an automorphism, a contradiction. $\Box$

\smallskip

 In conclusion, we recall a result of \cite[Theorem 1.3]{ShYu} saying that 
if, for a mapping $\varphi$  of ${\mathbb C}[x,y]$ with invertible Jacobian 
matrix, $\varphi(x)$ generates a proper retract of ${\mathbb C}[x,y]$, 
then $\varphi$ is an automorphism 
of  ${\mathbb C}[x,y]$. Then, the case where $\varphi(x)$ belongs to a
 proper retract but does not generate it, can be ruled out since in that 
case, $\varphi(x)=f(p(x,y)),$ where  $p(x,y)$ generates a proper retract of 
${\mathbb C}[x,y]$, and $f$ is some one-variable polynomial of degree $>$1. 
  The gradient of such a polynomial cannot form a row of any 
invertible Jacobian matrix, which can be easily seen from the ``chain rule" 
applied to $f(p(x,y))$. 

 Therefore, by Theorem 1 of the present paper, 
if $\varphi$ is a counterexample to the Jacobian conjecture for 
${\mathbb C}[x,y]$, then 
$\varphi(x)$ must be a test polynomial. Perhaps a way to prove the 
Jacobian conjecture for ${\mathbb C}[x,y]$ could be through showing that 
the  gradient  of a test polynomial cannot form a row of any 
invertible Jacobian matrix. This is known to be the case with (non-commutative) 
partial derivatives of a test element of a free group of rank 2, see 
\cite[Corollary 2.2.8]{ourbook}. \\

\noindent {\bf Acknowledgement}
\medskip 

  The first author is grateful to the  Department of Mathematics of the 
University of Hong Kong for its warm hospitality during his visit
when most of  this work was done. \\

\noindent 
 Department of Mathematics, The City  College  of New York, New York, 
NY 10031 
\smallskip

\noindent {\it e-mail address\/}: ~shpil@groups.sci.ccny.cuny.edu 

\smallskip

\noindent {\it http://www.sci.ccny.cuny.edu/\~{}shpil} \\

\noindent Department of Mathematics, The University of Hong Kong, 
Pokfulam Road, Hong Kong 

\smallskip

\noindent {\it e-mail address\/}: ~yujt@hkusua.hku.hk 
\smallskip

\noindent {\it http://hkumath.hku.hk/\~{}jtyu}

\end{document}